# Spaces of maps with symplectic graph

BY

Joseph Coffey

Courant Institute for Mathematical Sciences

*October 21, 2004*

**Abstract**

We consider the homotopy type of the space $\mathcal{M}^\sigma(\Sigma, \Gamma)$ of maps between symplectic surfaces $(\Sigma, \sigma_\Sigma)$ and $(\Gamma, \sigma_\Gamma)$ whose graphs form symplectic submanifolds of the product. We give a purely topological model for this space in terms of maps with constrained numbers of pre-images. We use this to show that the dependence of the homotopy type on the forms $\sigma_\Sigma$ and $\sigma_\Gamma$ is quantized- it changes only when the parameters pass certain discrete levels. When the domain $\Sigma$ is a sphere or torus, and further $\sigma_\Sigma(\Sigma) \leq \sigma_\Gamma(\Gamma)$, we compute the full homotopy type of the low degree components. We also give an example, showing that the homotopy type of the space of sections of a symplectic fibration $F$ must sometimes change as we deform $F$. Much of this work generalizes to $n$-dimensional manifolds equipped with volume forms.

## 1 Introduction and motivation

The spaces $\mathcal{H}ol(\mathbb{CP}^1, \mathbb{CP}^1)$ of holomorphic maps of $\mathbb{CP}^1$ to itself have enjoyed a long and fruitful study, beginning with Segal's work [Seg79]. One can view a holomorphic map in at least two ways: as a map preserving the complex structure, or as a map whose graph in $\mathbb{CP}^1 \times \mathbb{CP}^1$ is holomorphic. If one takes this second point of view we are led to consider maps whose graphs are otherwise constrained. In particular if one endows the domain $\Sigma$ and the range $\Gamma$ with volume forms $\sigma_\Sigma$ and $\sigma_\Gamma$, one can consider the maps $\mathcal{M}^\sigma(\Sigma, \Gamma)$ whose graphs in $\Sigma \times \Gamma$ are symplectic submanifolds with respect to the product symplectic structure:

$$\sigma = \pi_\Sigma^* \sigma_\Sigma + \pi_\Gamma^* \sigma_\Gamma$$

This paper is motivated by the study of homotopy type of $\mathcal{M}^\sigma(\Sigma, \Gamma)$. However the problem has a more natural setting in the product of two distinct $n$-manifolds, each equipped with a volume forms. Then $\mathcal{M}^\sigma(\Sigma, \Gamma)$ will be the space of maps whose graph $G$ is such that $\sigma|_G$ yields a volume form. In sections Spa81 and Spa81 we begin by by examining this more general problem.

In this context we construct a purely topological model of the homotopy type of $\mathcal{M}^\sigma(\Sigma, \Gamma)$ in terms of maps with constrained numbers of pre-images. We use this characterization to show that the dependence of the homotopy type of $\mathcal{M}^\sigma(\Sigma, \Gamma)$ on the form $\sigma$ is quantized; it may jump only at discrete intervals, described when the ratio $\frac{\sigma_\Sigma[\Sigma]}{\sigma_\Gamma[\Gamma]}$ passes an integer. This "quantized" topology is reminiscent of that found in the symplectomorphism groups of ruled surfaces. In fact this paper was originally motivated by efforts to understand these groups by primarily soft methods.

We also prove certain (general) identities between different components of $\mathcal{M}^\sigma(\Sigma, \Gamma)$ for different forms $\sigma$. Then in section 5 we use combine these identities with the theory of $J$-holomorphic spheres to show that the homotopy type of the space of sections of a fibration must sometimes change as we deform the fibration. This suggests both an interesting problem, and that the main results of this paper are probably out of reach of the methods of $J$-holomorphic curves. It also shows that any link with the homotopy type of symplectomorphism groups is complicated.

In section 4 we return to the case where both the domain and range are surfaces to compute the homotopy type of $\mathcal{M}^\sigma(\Sigma, \Gamma)$. We build concrete geometric models for these spaces, via a decomposition method which may be of independent interest. In an upcoming paper the we use these techniques to show the failure of the parametric $h$-principle for maps with prescribed Jacobian.

We succeed in providing a complete description of when the domain $\Sigma$ is a sphere or torus, and further $\sigma_\Sigma(\Sigma) \leq \sigma_\Gamma(\Gamma)$. Under these conditions, when $\Sigma$ is a sphere we see that the inclusion $\mathcal{H}ol(\mathbb{CP}^1,\Gamma) \hookrightarrow \mathcal{M}^\sigma(\mathbb{CP}^1, \Gamma)$ is a homotopy equivalence for maps of degree 0 or 1. However, when the domain is $T^2$ this fails dramatically. For instance the homology of the degree 0 component of $\mathcal{M}^\sigma(T^2, S^2)$ is infinitely generated in all dimensions.

## 1.1 Acknowledgements

This paper (save sections 4.5 and 5) was a portion of the authors PhD thesis, done under the direction of Dusa Mcduff and Dennis Sullivan. The author would also like to thank Pedro Santos for helpful discussion, the Instituto Superior Tecnico in Lisbon, the Courant Institute for Mathematical Sciences, and the Institute Hautes Etudes Scientifique for their hospitality. The author is supported in part by NSF grant DMS 0202553.

# 2 General Results

We consider the product of two $n$-manifolds: $\Sigma \times \Gamma$. We endow each factor with a volume form given by $\sigma_\Sigma$ and $\sigma_\Gamma$ respectively. These volume forms induce a product $n$-form on $\Sigma \times \Gamma$:

$$\sigma = \pi_\Sigma^* \sigma_\Sigma + \pi_\Gamma^* \sigma_\Gamma$$

We will call such an $n$-form on $\Sigma \times \Gamma$ a **product $n$-form**. Denote the space of $C^k$ maps from $\Gamma$ to $\Sigma$ by $\mathcal{M}(\Sigma, \Gamma)$. Those of degree $d$ we denote by $\mathcal{M}_d(\Sigma, \Gamma)$.

**Definition 2.1.** $f \in \mathcal{M}_d(\Sigma, \Gamma)$ *is a **positive map** if equivalently:*

1. *$\sigma$ restricted to the graph of $f$ is a volume form.*
2. *$f^* \sigma_\Gamma + \sigma_\Sigma$ is a volume form on $\Sigma$.*

We denote the space of positive maps $f \in \mathcal{M}_d(\Sigma, \Gamma)$ by $\mathcal{M}_d^\sigma(\Sigma, \Gamma)$. However when appropriate we will abbreviate; for example we may shorten $\mathcal{M}_d^\sigma(\Sigma, \Gamma)$ to $\mathcal{M}_d^\sigma$.

In this section and the next we will give a purely topological model for this space in terms of maps with constrained numbers of pre-images. For example we will show:

**Theorem. (Example of results)** *Consider the product of two $n$-manifolds: $\Sigma \times \Gamma$, with a product $n$-form $\sigma = \pi_\Sigma^* \sigma_\Sigma + \pi_\Gamma^* \sigma_\Gamma$. Suppose $\frac{\sigma_\Sigma[\Sigma]}{\sigma_\Gamma[\Gamma]} \leq 1$, then $\mathcal{M}_0^\sigma(\Sigma, \Gamma)$ is homotopy equivalent to the space of non-surjective maps from $\Sigma$ to $\Gamma$.*

In general, when the ratio $\frac{\sigma_\Sigma[\Sigma]}{\sigma_\Gamma[\Gamma]} > 1$ we will prove a similar result with the non-surjective maps replaced by "non Q-surjective" maps. We define these now in the next subsection.

## 2.1 A model for $\mathcal{M}_d^\sigma(\Sigma, \Gamma)$: Maps with constrained surjectivity

**Definition 2.2.** *A $C^1$ map $f \colon \Sigma \to \Gamma$ is called **non Q-Surjective** if there is an open ball in $U \subset \Gamma$ such that every $x \in U$ has less than $Q$ pre-images.*

Denote the space of $C^1$ smooth degree $d$, non $Q-$ surjective maps $\Sigma \to \Gamma$ by $\mathcal{NS}_d^Q$. Here $d \in \mathbb{Z}$, and $Q \in \mathbb{R}$. Note that $\mathcal{NS}_d^Q$ is empty if $Q \leq d$.

**Example 2.3.** $\mathcal{NS}_0^Q$ is the space of non-surjective maps for $0 < Q \leq 2$.

**Definition 2.4.** *For any product $n$-form $\sigma$ on $\Sigma \times \Gamma$, let $K(\sigma)$ denote the ratio:*

$$\frac{\sigma_\Sigma[\Sigma]}{\sigma_\Gamma[\Gamma]}$$

This section is devoted to the following purely topological description of $\mathcal{M}_\sigma(\Sigma, \Gamma)$ in in terms of the topology of $\Sigma$ and $\Gamma$ and this parameter $K(\sigma)$. When the form $\sigma$ is clear from context we will abbreviate $K(\sigma)$ by $K$. We also note that scaling the form $\sigma$ by a positive constant leaves $\mathcal{M}_\sigma(\Sigma, \Gamma)$ unchanged, thus we allow ourselves the following convention to simply the exposition:

**Convention:** *We scale each form $\sigma = \pi_\Sigma^* \sigma_\Sigma + \pi_\Gamma^* \sigma_\Gamma$ by a positive constant so that $\int_\Gamma \sigma_\Gamma = 1$. Then $K = \int_\Sigma \sigma_\Sigma$.*

**Theorem 2.5.** *The natural inclusion $i\colon \mathcal{M}_d^\sigma \hookrightarrow \mathcal{NS}_d^{2K(\sigma)+d}$ is a deformation retract.*

For $x \in \mathbb{R}$ we denote by $\lceil x \rceil$ the least integer $n \in \mathbb{Z}$ such that $n \geq x$. Then this characterization has the following corollaries:

**Corollary 2.6.** *( The dependence of $\mathcal{M}^\sigma$ on $\sigma$ is quantized) $\mathcal{M}_d^{\sigma_1}$ is homotopy equivalent to $\mathcal{M}_d^{\sigma_2}$ if $\lceil K(\sigma_1) \rceil = \lceil K(\sigma_2) \rceil$.*

**Proof.** The corollary follows from the observation that $\mathcal{NS}_d^Q$'s dependence on $Q$ is quantized. More precisely we observe that $\mathcal{NS}_d^{Q_1} = \mathcal{NS}_d^{Q_2}$ if
$$\left\lceil \frac{Q_1 - |d|}{2} \right\rceil = \left\lceil \frac{Q_2 - |d|}{2} \right\rceil$$

To see this let $f \in \mathcal{NS}_d^Q$. Let $U \subset \Gamma$ be an open ball such that every $x \in U$ has less than $Q$ pre-images. Let $x \in U$ be a regular value of $f$. As $\Gamma$ and $\Sigma$ are both compact, the set of regular values of $f$ is open. Let $U_{reg}$ be a neighborhood of $x$ consisting of regular values of $f$. Then for every $y$ in $U$, $f^{-1}(y)$ has the same cardinality.
$$\mathrm{Card}(f^{-1}(y)) = |d| + 2l < Q$$
for $l$ a positive integer. The stated quantization of $\mathcal{NS}_d^Q$ then follows, and thus the corollary. For if
$$\left\lceil \frac{2K(\sigma_1) + d - |d|}{2} \right\rceil = \left\lceil \frac{2K(\sigma_2) + d - |d|}{2} \right\rceil$$
then $\lceil K(\sigma_1) \rceil = \lceil K(\sigma_2) \rceil$ as $\frac{d-|d|}{2} \in \mathbb{Z}$. □

**Example 2.7.** Consider $\mathcal{NS}_1^{2K+1}$. A regular value of a degree one map can have any odd number of pre-images: 1,3,5,7.... Thus $\mathcal{NS}_1^{2K+1}$ is empty when $K < 0$, and changes only when $K$ passes an integer.

**Corollary 2.8.** *Suppose that either the domain $\Sigma$ or range $\Gamma$ admits a degree $-1$ diffeomorphism $\gamma$, then if $K(\sigma_2) = K(\sigma_1) + d$, $\mathcal{M}_d^{\sigma_1}$ is homotopy equivalent to $\mathcal{M}_{-d}^{\sigma_2}$.*

**Proof.** By theorem 2.5 we have homotopy equivalences:
$$\mathcal{M}_d^{\sigma_1} \simeq \mathcal{NS}_d^{2K(\sigma_1)+d}$$
$$\mathcal{M}_{-d}^{\sigma_2} \simeq \mathcal{NS}_{-d}^{2K(\sigma_2)-d} = \mathcal{NS}_{-d}^{2K(\sigma_1)+d}$$

Suppose that the range $\Gamma$ possesses the orientation reversing diffeomorphism. Then $\phi$ determines a homeomorphism:
$$\gamma_*\colon \mathcal{NS}_{-d}^{2K(\sigma_1)+d} \to \mathcal{NS}_d^{2K(\sigma_1)+d}$$
$$f \to \gamma \circ f$$

yielding the corollary. If $\phi$ is defined on $\Sigma$ we can define $\gamma_*$ by precomposition. □

**Example 2.9.** Suppose $\sigma_1$ and $\sigma_2$ are product 2-forms on $S^2 \times S^2$, such that $K(\sigma_2) = K(\sigma_1) + 1$. Then, as $S^2$ admits an orientation reversing involution $\gamma$, $\mathcal{M}_1^{\sigma_1}(S^2, S^2)$ is homotopy equivalent to $\mathcal{M}_{-1}^{\sigma_2}(S^2, S^2)$.

# 3 Proof of Theorem 2.5

## 3.1 Strategy for construction of $i\colon \mathcal{M}_d^\sigma \hookrightarrow \mathcal{NS}_d^{2K(\sigma)+d}$

We will seek to use Moser's Lemma to the greatest extent possible. First we will apply it "in the domain" to show that $\mathcal{M}_d^\sigma$ is a deformation retract of $\mathcal{NV}_d^{K(\sigma)}$ – those $f \in \mathcal{M}_d(\Sigma, \Gamma)$ which have bounds on how great an area of the domain $f$ maps with negative orientation (subsections 3.2 - 3.3.1). Then we will apply Moser's lemma "in the range" to show that $\mathcal{NV}_d^{K(\sigma)}$ is a deformation retract of $\mathcal{NS}_d^{2K(\sigma)+d}$.

## 3.2 Maps of bounded negative volume

### 3.2.1 Definition of negative volume

**Definition 3.1.** *Let $\alpha$ be an n-form on an oriented n-manifold $\Sigma$.*
$$\mathcal{N}(\alpha) = \{x \in \Sigma \colon \alpha(x) < 0\}$$
*where this sign is determined by the orientation of M.*

If $f\colon \Sigma \to \Gamma$ and $\sigma_\Gamma$ is a volume form on on $\Gamma$ we will sometimes denote $\mathcal{N}(f^*\sigma_\Gamma)$ by $\mathcal{N}(f)$, for this set depends only on $f$, and not on the choice of volume form $\sigma_\Gamma$. $\mathcal{N}(f)$ is the subset of the domain $\Sigma$ where $f$ reverses orientation.

**Definition 3.2.** *Define the **negative volume** of an n-form $\alpha$ denoted $NV(\alpha)$ to be $-\int_{\mathcal{N}(\alpha)} \alpha$.*

If $f\colon \Sigma \to \Gamma$ we will denote $NV(f^*\sigma_\Gamma)$ by $NV(f)$.

**Definition 3.3.** *Denote by $\mathcal{NV}_d^K$ those $f \in \mathcal{M}_d$ such that $NV(f^*\sigma_\Gamma) < K$. We will refer to these as **maps of bounded negative volume**.*

### 3.2.2 Basic properties of negative volume

For a map $f \in \mathcal{M}(\Sigma, \Gamma)$ denote by $reg(f) \subset \Gamma$ the regular values of $f$. For $x \in reg(f)$ denote by $\mu_f(x)$ the cardinality of the set $f^{-1}(x) \cap \mathcal{N}(f)$. Then:

**Lemma 3.4.** *Let $\eta$ be a volume form on $\Gamma$. Then $NV(f^*\eta) = \int_{reg(f)} \mu_f(x) \eta$*

**Proof.** $f|_{\mathcal{N}(f)}$ is a covering map over each connected component $X_i$ of $reg(f)$, the set of regular values of $f$ in $X$. This may be the empty cover over certain components $X_i$ - some regular $x$ may have no negative pre-images. $\mu_f(x)$ is constant for $x \in X_i$, and gives the number of sheets in this cover. Thus
$$\int_{\mathcal{N}(f) \cap f^{-1}(X_i)} f^*\eta = \mu_f(x) \int_{X_i} \mu$$
We gain the Lemma by integrating over each $X_i$. □

**Lemma 3.5.** *(Invariance of domain) Let $\gamma \in \mathrm{Diff}(\Sigma)$, $\alpha$ be an n-form on $\Sigma$. Then $NV(\gamma^*\alpha) = NV(\alpha)$*

**Proof.** $\mathcal{N}(\gamma^*\alpha) = \gamma^{-1}(\mathcal{N}(\alpha))$. Thus
$$\int_{\mathcal{N}(\gamma^*\alpha)} \gamma^*\alpha = \int_{\gamma^{-1}(\mathcal{N}(\alpha))} \gamma^*\alpha = \int_{\mathcal{N}(\alpha)} \alpha \qquad \square$$

**Lemma 3.6.** *(Continuity) $NV\colon \Omega^n(\Sigma) \to \mathbb{R}$ is continuous in the $C^0$ topology on forms.*

**Proof.** This follows immediately from the observation that:
$$NV(\alpha) = \frac{1}{2} \int_\Sigma (|\alpha| - \alpha) \qquad \square$$

### 3.3 Positive maps are a deformation retract of maps with bounded negative volume

This section is devoted to the proof of the following proposition:

**Proposition 3.7.** *The natural inclusion* $j\colon \mathcal{M}_d^\sigma \hookrightarrow \mathcal{NV}_d^{K(\sigma)}$ *is a deformation retract.*

**Proof.**

Note that $\mathcal{M}_d^\sigma \subset \mathcal{NV}_d^{K(\sigma)}$. For if $f \in \mathcal{M}_d^\sigma$ the following equation holds for any (measurable) $U \subset \Sigma$:

$$\int_U (\sigma_\Sigma + f^*\sigma_\Gamma) > 0$$

If we take $U = \mathcal{N}(f)$ we have:

$$\left(\int_{\mathcal{N}(f)} \sigma_\Sigma\right) - NV(f) > 0$$

and thus a bound on the negative volume of $f$:

$$NV(f) < \int_{\mathcal{N}(f)} \sigma_\Sigma < \int_\Sigma \sigma_\Sigma = K$$

We now set upon proving that this inclusion is a weak deformation retract. Consider a disc $\rho$ of non-surjective maps with boundary in $\mathcal{M}_d^\sigma$:

$$\rho\colon (D^n, \partial D^n) \to (\mathcal{NV}_d^K, \mathcal{M}_\sigma^d)$$

We will construct a retraction of this disc into the space of positive maps $\mathcal{M}_\sigma^d$. i.e. we will construct a homotopy of pairs

$$\rho_t\colon (D^n, \partial D^n) \to (\mathcal{NV}_d^K, \mathcal{M}_\sigma^d) \quad t \in [0,1]$$

such that:

1. $\rho_0(x) = \rho$ for each $x \in D^n$.
2. $\rho_1(D^n) \subset \mathcal{M}_\sigma^d$.
3. $\rho_t|_{\partial D^n} = \rho|_{\partial D^n}$ for all $t$.

Denote the volume forms on $\Sigma$ in class $[\sigma_\Sigma]$ by $\mathcal{V}_K$. We will construct $\rho_t$ by constructing a family of volume forms:

$$\tau_\rho\colon (D^m, \partial D^m) \to (\mathcal{V}_K, \sigma_\Sigma)$$

such that each map $\rho(x)$ is positve with respect to the product $n$-form $\sigma_\tau(x) = \pi_\Sigma^*(\tau_\rho(x)) + \pi_\Gamma^*(\sigma_\Gamma)$:

$$\rho(x) \in \mathcal{M}_d^{\sigma_\tau(x)}$$

In 3.3.2 we will discuss how Moser's Lemma will provide a family of diffeomorphisms $M_\tau(x, 1)$ of $\Sigma$ such that $M_\tau(x, 1)^*(\sigma_\Sigma) = \tau_\rho(x)$. Our homotopy is then given by:

$$\rho_t(x) = \rho(x) \circ M_\tau^{-1}(x, t)$$

and it provides the necessary retraction of $\mathcal{NV}_d^K$ into $\mathcal{M}_\sigma^d$.

#### 3.3.1 Construction of family of forms on the domain $\Sigma$

**Lemma 3.8.** *Let $\rho\colon D^m \to \mathcal{NV}_d^{K(\sigma)}$ be a disc of maps, such that $\rho(\partial D^m) \subset \mathcal{M}_\sigma^d$. Then there is a continuous map $\tau_\rho\colon D^m \to \mathcal{V}_K$ such that:*

1. $\rho(x)^*\sigma_\Gamma + \tau_\rho(x) > 0$

2. $\tau_\rho(\partial D^m) = \sigma_\Sigma$

**Proof.**

**Step 1: Constructing $\tau_\rho$ for a single fixed map** We begin by fixing a $g \in \mathcal{NV}_d^K$, and construct a $\tau_g \in \mathcal{V}_K$ such that:

$$g^*\sigma_\Gamma + \tau_g > 0$$

We build the form $\tau_g$ from $\sigma_\Sigma$ by first modifying $\sigma_\Sigma$ within $\mathcal{N}(g)$ to attain the above inequality, then we scale $\sigma_\Sigma$ away from $\mathcal{N}(g)$ so that the resulting form has the proper cohomology class, and thus lies in $\mathcal{V}_K$. There is enough volume to go around precisely because $NV(g) < K(\sigma)$. To wit:

Let $v_\delta$ be a volume form on $\Sigma$ such that:

$$\int_\Sigma v_\delta < K(\sigma) - NV(g)$$

Let $U_\delta$ be an open neighborhood of $\mathcal{N}(g)$ such that $v_\delta - g^*\sigma_\Gamma > 0$, and let $\psi$ be a $C^\infty$ function such that:

1. $\psi(x) = 0$ for $x \in \Sigma \backslash U_\delta$
2. $\psi = 1$ for $x \in \mathcal{N}(g)$

Then for $\tau_g' = \psi \cdot (v_\delta - g^*\sigma_\Gamma) + (1 - \psi) \cdot \sigma_\Sigma$ we achieve:

$$g^*\sigma_\Gamma + \tau_g' > 0$$

Moreover we note that

$$\int_{\mathcal{N}(g)} \tau_g' = \int_{\mathcal{N}(g)} (v_\delta - g^*\sigma_\Gamma) = \left(\int_{\mathcal{N}(g)} v_\delta\right) + NV(g) < K(\sigma)$$

if we choose $v_\delta$ to be small. We now scale $\tau_g'$ in the complement of $\mathcal{N}(g)$ so so that the resulting forms are in $\Omega_\sigma$. We now define $\varphi$ as $C^\infty$ function on $\Sigma$ such that:

1. $\varphi(x) = 1$ if $x \in \overline{\mathcal{N}(g)}$
2. $\varphi(d, x) < 1$ elsewhere.

For $s \in [1, \infty)$ we examine the integral:

$$F(k) = \int_\Sigma (\varphi^s \cdot \tau_g')$$

This a monotone, decreasing, continuous function of the real parameter s. As $s$ approaches $\infty$, $\psi^s$ converges to the characteristic function of $\overline{\mathcal{N}(g)}$. Thus:

$$\lim_{s \to \infty} F(s) = \int_{\mathcal{N}(g)} \tau_g' < K(\sigma) = \int_\Sigma \sigma_\Sigma$$

Thus for some value $s = t$, $F(t) = K(\sigma)$, and $(\varphi^t \cdot \tau_g') \in \Omega_\sigma$. Since our scaling took place entirely away from $\mathcal{N}(g)$ we retain the inequality:

$$g^*\sigma_\Gamma + \varphi^t \cdot \tau_g' > 0$$

and $\tau_g = \varphi^t \cdot \tau_g'$ is the required form.

**Step 2: Constructing $\zeta_\rho$ in families via partition of unity:** We now globalize the previous construction to the disc of maps $\rho$ via a partition of unity.

Let

$$\bigcup_{i \in I} U_i$$

be a finite covering of $D^n$ by open sets $U_i$ with forms $\tau_{\rho(x_i)}$ such that for $x \in U_i$ $\rho(x) \in \mathcal{M}_d^{\sigma_\tau(x_i)}$, where

$$\sigma_\tau(x) = \pi_\Sigma^*(\tau_{\rho(x_i)}) + \pi_\Gamma^* \sigma_\Gamma$$

Let $U_{\partial D}$ be a neighborhood of the boundary such that for $x \in U_{\partial D}$ $(x) \in \mathcal{M}_d^\sigma$. Finally let $\{\phi_i, \phi_\partial\}$ be a partition of unity subordinate to the covering of $D^n$ given by the $U_i$ and $U_{\partial D}$. Then:

$$\tau_\rho(x) = \sum_{i \in I} \phi_i(x) \cdot \tau_{\rho(x_i)} + \phi_\partial \sigma_\Sigma$$

is the family of forms required. The proof of Lemma 3.8 is thus completed. □

### 3.3.2 Applying Moser's Lemma on the domain Σ

$\mathcal{V}_K$ forms a convex set. Thus we can find a homotopy $\tau_t$ of $\tau_\rho$ to the constant sphere, based at $\sigma_\Sigma$ within $\mathcal{V}_K$. Moser's Lemma applies and so if we denote the diffeomorphisms of $\Sigma$ by $Diff(\Sigma)$ we obtain:

$$M_\tau: (D^n \times I, \partial D \times I) \to (Diff(\Sigma), Id)$$

such that $M_\tau(x,1)^*(\sigma_\Sigma) = \tau_\rho(x)$. Let

$$\rho_t(x) = \rho(x) \circ M_\tau^{-1}(x,t)$$

Then

$$NV(\rho_t(x)) = NV\left(\rho(x) \circ M_\tau^{-1}(x,t)\right) = NV(\rho(x)) < K$$

by Lemma 3.5 (Invariance of domain). Further we see that $\rho_t$ gives a retraction into $\mathcal{M}_d^\sigma$, for

$$\rho_1(x)^*(\sigma_\Gamma) + \sigma_\Sigma = M_\tau^{-1}(x,1)^* \rho(x)^*(\sigma_\Gamma) + \sigma_\Sigma$$

is indeed a volume form by pulling it back by $M_\tau(x,1)$:

$$M_\tau(x,1)^*(\rho_1(x)^*(\sigma_\Gamma) + \sigma_\Sigma) = \rho(x)^*(\sigma_\Gamma) + M_\tau(x,1)^* \sigma_\Sigma$$
$$= \rho(x)^*(\sigma_\Gamma) + \tau_\rho(x)$$

which is a volume form by condition 1 in Lemma 3.8. □

## 3.4 Maps of bounded negative volume are a weak deformation retract of non $Q$-surjective maps

For $f \in C^1(\Sigma, \Gamma)$ denote by $\mu_f(x)$ the cardinality of $f^{-1}(x) \cap \mathcal{N}(f)$.

**Proposition 3.9.** *The natural inclusion $i: \mathcal{NV}_d^K \hookrightarrow \mathcal{NS}_d^{2K+d}$ is a weak deformation retract.*

We first show that $\mathcal{NV}_d^K \subset \mathcal{NS}_d^{2K+d}$. Let $f$ be a degree $d$ map in $C^1(\Sigma, \Gamma) \backslash \mathcal{NS}_d^{2K+d}$. Then I claim that for every regular value $x \in reg(f)$, $\mu_f(x) > K$: For if $d \geq 0$ then $x$ has (at least) $2K$ "excess" pre-images. Half of these must be negative. If $d < 0$ then $x$ has (at least) $2K - 2|d|$ "excess" pre-images. Again half of these must be negative. We also have $|d|$ negative pre-images coming from the degree of the map. This again yields $K = K - |d| + |d|$ in total.

By Lemma 3.4,

$$\begin{aligned} NV(f) &= \int_{reg(f)} \mu_f \sigma \\ &= \int_{reg(f) \cap X} \mu_f \sigma \\ &\geq K \int_{reg(f) \cap X} \sigma \\ &= K \end{aligned}$$

And so $\mathcal{NV}_d^K \subset \mathcal{NS}_d^{2K+d}$.

We now perform the retraction. Consider a disc:
$$\rho \colon (D^n, \partial D^n) \to (\mathcal{NS}_d^{2K+d}, \mathcal{NV}_d^K)$$

We will construct a retraction of this disc into the space of positive maps $\mathcal{M}_d^\sigma$. i.e. we will construct a homotopy of pairs

$$\rho_t \colon (D^n, \partial D^n) \to (\mathcal{NV}_d^K, \mathcal{M}_\sigma^d) \ t \in [0,1]$$

Our strategy is same as before. We construct a family of forms and then apply Moser's Lemma. We remind the reader that for an $n$-form $\alpha$ we define $NV(\alpha)$ to be $-\int_{\mathcal{N}(\alpha)} \alpha$.

### 3.4.1 Constructing a family of form on the range $\Gamma$

**Lemma 3.10.** *Let $\rho \colon (D^n, \partial D) \to (\mathcal{NS}_d^{2K+d}, \mathcal{NV}_d^K)$ There is a continuous map of pairs $\tau_\rho \colon (D^n, \partial D) \to (\Omega_\Gamma, \sigma_\Gamma)$ such that:*

$$NV(\rho(x)^* \tau_\rho(x)) < K$$

The proof of this Lemma proceeds along the same lines as that of Lemma 3.8.

**Step 1: Constructing $\tau_\rho$ for a single fixed map $g$** Partition the range $\Gamma$ into a set $\Gamma_-$ with less than $2K + d$ pre images, and its complement $\Gamma_+$. Partition the domain into $\Sigma_+ = g^{-1}(\Gamma_+)$ and $\Sigma_- = g^{-1}(\Gamma_-)$. Finally, denote the volume forms $\tau$ on $\Gamma$ such that $\int_\Gamma \tau = \int_\Gamma \sigma_\Gamma = 1$ by $\mathcal{V}$.

As $g \in \mathcal{NS}_d^{2K+d}$, $\Gamma_-$ has nonempty interior we may thus find a volume form $\tau_g^\epsilon \in \mathcal{V}$ such that $\Gamma_-$ very large:

$$\int_{\Gamma_-} \tau_g^\epsilon = 1 - \delta$$

while

$$\frac{\tau_g^\epsilon(x)}{\sigma_\Gamma(x)} < \epsilon$$

for all $x \in \Gamma_+$ where $\epsilon > 0$ is a constant which can be made as small as we like by making $\delta$ small.

Now $x \in \Gamma_-$ has $< 2K + d$ pre images under a map and thus the number of negative pre-images $\mu_g(x)$ is a also bounded:

$$\mu_g(x) < K - \epsilon_1$$

for some $\epsilon_1 > 0$. By lemma 3.4:

$$\begin{aligned} NV(g^* \tau_g^\epsilon |_{\Sigma_-}) &= \int_{\Gamma_-} \mu_g(x) \tau_g^\epsilon \\ &< (K - \epsilon_1) \int_{\Gamma_-} \tau_g^\epsilon \\ &< (K - \epsilon_1)(1 - \delta) \end{aligned}$$

On $\Sigma_+$ we have the following bound:

$$NV(g^* \tau_g^\epsilon |_{\Sigma_+}) < \epsilon \Big( NV(g^* \sigma_\Gamma |_{\Sigma_+}) \Big)$$

As:

$$NV(g^* \tau_g^\epsilon) = NV(g^* \tau_g^\epsilon |_{\Sigma_-}) + NV(g^* \tau_g^\epsilon |_{\Sigma_+})$$

We have:

$$NV(g^* \tau_g^\epsilon) < (K - \epsilon_1)(1 - \delta) + \epsilon(NV(g^* \sigma_\Gamma |_{\Sigma_+}))$$

which approaches $K - \epsilon_1$ as $\epsilon, \delta \to 0$. We can thus make $NV(g^* \tau_g^\epsilon) < K$ by making both $\delta$ and $\epsilon$ small.

**Step 2: Constructing $\tau_\rho$ in families via partition of unity**

By lemma 3.6, for a fixed $x_i \in D^n$, $NV(g^*\tau_i): C^1(\Sigma, \Gamma) \to \mathbb{R}$ is a continuous function in $g$. Thus, after applying the previous step for each map $\rho(x)$, one can find a covering

$$\bigcup_{i \in I} U_i$$

$D^n$ by open sets $U_i$ with forms $\tau_i$ such that for $x \in U_{i_i}$ $NV(\rho(x)^*\tau_i) < K$. By compactness of $D^n$ we can ensure that $I$ is finite by taking a finite subcover.

Let $U_{\partial D}$ be a neighborhood of the boundary such that $NV(\rho(x)^*\sigma_\Gamma) < K$ for $d \in U_{\partial D}$. Finally let $\{\phi_i, \phi_\partial\}$ be a partition of unity subordinate to the covering of $D^n$ given by the $U_{d_i}$ and $U_{\partial D}$. Then, as both $NV(f^* \cdot) < K$ and having cohomology class $[\sigma_\Gamma]$ are convex conditions on $n$-forms:

$$\tau_\rho(x) = v_\partial(x)\sigma_\Gamma + \sum_{i \in I} \phi_i(x)\tau_i$$

satisfies the conditions of the proposition.

### 3.4.2 Applying Moser's Lemma on the range $\Gamma$

We can find a homotopy $\tau_t$ of $\tau_\rho$ to the constant sphere mapping all of $D^n$ to $\sigma_\Gamma$ within volume forms on $\Gamma$ with cohomology class $[\sigma_\Gamma]$ Moser's Lemma yields a family of diffeomorphisms:

$$M_\tau: (D^n \times I, \partial D \times I) \to (Diff(\Sigma), Id)$$

such that $M_\tau(x, 1)^*(\sigma_\Gamma) = \tau_\rho(x)$. Let

$$\rho_t(x) = M_\tau(x, t) \circ \rho(x)$$

Clearly $\rho_t$ remains in $\mathcal{NS}_d^K$: post-composing a map with a diffeomorphism does not change its $Q$-surjectivity.

I claim that

$$NV(\rho_1(x)^*\sigma_\Gamma) < K$$

For $\rho_1(x) = M_\tau(x, 1) \circ \rho(x)$. Thus $\rho_1(x)^*\sigma_\Gamma = \rho(d)^*M_\tau(x, 1)^*\sigma_\Gamma$. So

$$\begin{aligned} NV(\rho_1(x)^*\sigma_\Gamma) &= NV(\rho(x)^*M_\tau(x,1)^*\sigma_\Gamma) \\ &= NV(\rho(x)^*\tau_\rho(x)) < K \end{aligned}$$

Combining Propositions 3.7 and 3.9 we achieve Theorem 2.5.

## 4 Computations- Making geometric models when the domain $\Sigma$ and range $\Gamma$ are surfaces

In this section we examine the case where $\Sigma$ and $\Gamma$ are both symplectic surfaces. Then $\sigma = \pi_\Sigma^*\sigma_\Sigma + \pi_\Gamma^*\sigma_\Gamma$ is a symplectic form, and the positive sections $\mathcal{M}_\sigma^d$ are symplectic sections for the product fibration. In certain situations we can completely describe the homotopy type of $\mathcal{M}_\sigma^d(\Sigma, \Gamma)$.

### 4.1 Description of Results

We now describe the results of our computations. The proofs of the main result we postpone to the next subsections, however we do derive a few corollaries.

#### 4.1.1 Case: Domain is a sphere

**Theorem 4.1.** *Let $(\Gamma, \sigma_\Gamma)$ and $(S^2, \sigma_\Sigma)$ be symplectic surfaces. Then $\mathcal{M}_0^\sigma(S^2, \Gamma)$ is homotopy equivalent to $\Gamma$ for $K(\sigma) \in (0, 1)$*

**Corollary 4.2.** *The homotopy type of $\mathcal{M}_0^\sigma(S^2, S^2)$ is a non-constant function of $K(\sigma)$.*

**Theorem 4.3.** *Let $(\Sigma = S^2, \sigma_\Sigma)$ and $(\Gamma = S^2, \sigma_\Gamma)$ be symplectic surfaces. $\mathcal{M}_d^\sigma(S^2, S^2)$ is homotopy equivalent to $SO(3)$ for:*

1. *$K(\sigma) \in (0, 1]$, $d = 1$*
2. *$K(\sigma) \in (1, 2]$, $d = 1$*

The proofs of these Theorems are given in subsection 4.4.

### 4.1.2 Case: Domain is a Torus

We now examine the case when the domain $\Sigma$ is the torus $T^2$. We begin by fixing some initial data. We will use this data to describe our model for $\mathcal{M}_0^\sigma(T^2, \Gamma)$:

**Definition 4.4.** *Denote by $\hat{\mathbb{Q}}$ the set of ordered pairs $(a, b)$ of relatively prime non zero integers, along with the pairs $(0, 1)$ and $(1, 0)$.*

**Definition 4.5.** *Let $\gamma_1$ and $\gamma_2$ denote a pair of simple closed curves in $T^2$, such that $\gamma_1 \cap \gamma_2 = x \in T^2$, and such that $T^2 \backslash \gamma_1 \cap \gamma_2$ is a disk. Let $z \in S^1$. For each pair $(a, b) \in \hat{\mathbb{Q}}$ let $\phi_{a,b}$ be a map $T^2$ to $S^1$ such that $\phi_{a,b}|_{\gamma_1}$ has degree $a$, and $\phi_{a,b}|_{\gamma_2}$ has degree $b$, and $\phi_{a,b}(x) = z$. Further we choose these retractions so that the maps $\phi_{0,1}(\gamma_1) = z$ and $\phi_{1,0}(\gamma) = z$.*

**Definition 4.6.** *Let $\Theta$ be a (possibly non compact) surface. Denote by $\Omega_\Theta$ the maps $f: T^2 \to \Theta$ such that $f$ factors as $f = g \circ \phi_{a,b}$ for some $a, b \in \mathbb{Z}$, $g \in \mathcal{M}(S^1, \Theta)$.*

Note that the only maps which admit more than one such factorization are the constant maps.

**Theorem 4.7.** *The inclusion $j: \Omega_\Gamma \hookrightarrow \mathcal{N} \mathcal{S}_0^1(T^2, \Gamma)$ is a homotopy equivalence. Thus, if $K(\sigma) \leq 1$, $\mathcal{M}_0^\sigma(T^2, \Gamma)$ is homotopy equivalent to $\Omega_\Gamma$.*

The proof of this Theorem is in subsection 4.5. We note that it implies that $\mathcal{M}_0^\sigma(T^2, \Gamma)$ is homotopy equivalent to $\Omega_\Gamma$ and we establish the following corollaries. The first is just a more geometric description of $\Omega_\Gamma$. It is immediate.

**Corollary 4.8.** *Let $(T^2, \sigma_\Gamma)$ and $(S^2, \sigma_\Gamma)$ be symplectic surfaces such that $K(\sigma) \leq 1$. Then $\mathcal{M}_0^\sigma(T^2, S^2)$ is homotopy equivalent to countably many copies of the free loop space of $S^2$, joined over the constant loops:*

$$\coprod_{i \in \hat{\mathbb{Q}}} \mathcal{M}(S_i^1, S^2)\} / \sim$$

*where two maps $f_i \in \mathcal{M}(S_i^1, S^2)$ and $f_j \in \mathcal{M}(S_j^1, S^2)$ are equivalent if there is a $z \in S^2$ such that $f_i(S_i^1) = f_j(S_j^1) = z$.*

We can use this description to compute both $\pi_1(\mathcal{M}_0^\sigma(T^2, S^2))$ and $H_*(\mathcal{M}_0^\sigma(T^2, S^2))$:

**Corollary 4.9.** *Let $(T^2, \sigma_\Sigma)$ and $(S^2, \sigma_\Gamma)$ be such that $K(\sigma) \leq 1$ then $\mathcal{M}_0^\sigma(T^2, S^2)$ has infinitely generated $\pi_1$, and infinite dimensional homology in each dimension.*

**Proof.** These are immediate consequences of van Kampen's theorem and Mayer Vietoris respectively. □

**Corollary 4.10.** *Let $(T^2, \sigma_\Sigma)$ and $(\Gamma, \sigma_\Gamma)$ be symplectic surfaces such that $K(\sigma) \leq 1$ then*

1. *$\mathcal{M}_0^\sigma(T^2, S^2) \hookrightarrow \mathcal{M}_0(T^2, S^2)$ is not a homotopy equivalence, and thus $\mathcal{M}_0^\sigma(T^2, S^2)$ is a nontrivial function of $K(\sigma)$.*
2. *If $\Gamma$ is a surface of genus $> 0$, $\mathcal{M}_0^\sigma(T^2, \Gamma) \hookrightarrow \mathcal{M}_0(T^2, \Gamma)$ is a homotopy equivalence.*

**Proof.**

**Case: range is $S^2$.** The inclusion $\mathcal{M}_0^\sigma(T^2, S^2) \hookrightarrow \mathcal{M}_0(T^2, S^2)$ cannot by a homotopy equivalence as $\mathcal{M}^0(T^2, S^2)$ has finitely generated $\pi_1$, while $\mathcal{M}_0^\sigma(T^2, \Gamma)$ does not by Corollary 4.9.

**Case: range is higher genus.** By Theorem 4.7 the inclusion $j\colon \Omega_\Gamma \hookrightarrow \mathcal{NS}_0^1(T^2, \Gamma)$ is a homotopy equivalence. Thus is it sufficient to show that $i\colon \Omega_\Gamma \hookrightarrow \mathcal{M}_0(T^2, \Gamma)$ is a homotopy equivalence. If we fiber both $\Omega_\Gamma$ and $\mathcal{M}_0(T^2, \Gamma)$ over $\Gamma$ by $f \to f(x)$ we have the following morphism of fibrations.

$$\begin{array}{ccc} (\Omega_\Gamma, y) & \xrightarrow{i_{\text{fb}}} & (\mathcal{M}_0, y) \\ \downarrow & & \downarrow \\ \Omega_\Gamma & \xrightarrow{i} & \mathcal{M}_0 \\ \downarrow & & \downarrow \\ \Theta & \xrightarrow{\text{id}} & \Theta \end{array}$$

(We have suppressed the $T^2$ and $\Gamma$ from our notations.) Both spaces are connected thus it is sufficient to show that $i_{\text{fb}}$ is a homotopy equivalence. The result will then follow by the five lemma [Spa81].

Thus the corollary follows from the following Lemma, whose proof we postpone to 4.5.3. It is also the key fact underlying Theorem 4.7.

**Lemma 4.11.** *Let $\Theta$ be a (potentially non compact) surface. If $\Theta \neq S^2$, then the inclusion $j\colon (\Omega_\Theta, y) \hookrightarrow (\mathcal{M}^0(T^2, \Theta), y)$ is a homotopy equivalence.*

$\square$

## 4.2 Strategy of proof- Chopping up spaces

Our main tool is the following elementary (but not so well known) lemma in homotopy theory:

**Proposition 4.12.** *(Homotopy Decomposition Lemma[Gra75]) Let $f\colon X \to Y$ be a continuous map. Let $\bigcup_{j \in J} U_Y^j$ be a finite covering of $Y$ by open sets, and denote $f^{-1}(U_Y^i)$ by $U_X^i$. Suppose that for each $J \subset I$ the the restriction*

$$f\colon \bigcap_{j \in J} U_X^j \to \bigcap_{j \in J} U_Y^j$$

*is a homotopy equivalence. Then $f$ is a homotopy equivalence*

We calculate the homotopy type of each space $Y$ by

1. Finding an appropriate model space $X \subset Y$ (like the constant maps, or $\Omega_\Gamma$. $X$ should be simpler, and have an understandable homotopy type.)

2. Finding an appropriate covering of $Y$ so that it satisfies the hypotheses of the Homotopy Decomposition Lemma.

One can think of this as generalization of a CW decomposition of a space. We divide each space into pieces that map be complicated, but the map $f$ gives an equivalence of each piece and all their intersections.

## 4.3 Discrete Approximation

In most cases it is necessary to replaces the space $Y$ with a sequence of discrete approximations. We describe these now.

Let $T_0$ be a triangulation of the range $\Gamma$. We consider the system of triangulations $T_i$ of $\Gamma$, where $T_{i+1}$ denotes the barycentric subdivision of $T_i$. Endow $\Gamma$ with an auxiliary Riemannian metric $g$.

**Definition 4.13.** *(Nets of balls) Let $\mathcal{B}_i$ denote the set of balls $B(v_{ij}, \epsilon_i)$ centered at vertices $v_{ij}$ of $T_i$ and with radius $\epsilon_i > 0$ such that:*

1. $\epsilon_{i+1} < \epsilon_i$
2. *The balls $B(v_{ij}, \epsilon_i)$ are disjoint.*

**Definition 4.14.** *Denote by $(\mathcal{NS}_d^Q, \mathcal{B}_i)$ the degree $d$ maps $f \in C^1(\Sigma, \Gamma)$ such that there is a ball $\Delta \in \mathcal{B}_i$ such that $f^{-1}(\Delta)$ consists of $k < Q$ disjoint balls $D_i$, and $f\colon D_i \to \Delta$ is a diffeomorphism, and $f$ maps $\Sigma \setminus \bigcup D_i \to \Gamma \setminus \Delta$.*

*If $Q = 0$ these are the degree $0$ maps which miss $\Delta$.*

**Lemma 4.15.** *(Discrete approximation)* $\mathcal{NS}_d^Q$ *is the direct limit of the system*

$$(\mathcal{NS}_d^Q, \mathcal{B}_0) \hookrightarrow (\mathcal{NS}_d^Q, \mathcal{B}_1) \hookrightarrow \cdots (\mathcal{NS}_d^Q, \mathcal{B}_i) \hookrightarrow (\mathcal{NS}_d^Q, \mathcal{B}_{i+1}) \hookrightarrow \cdots \hookrightarrow \mathcal{NS}_d^Q$$

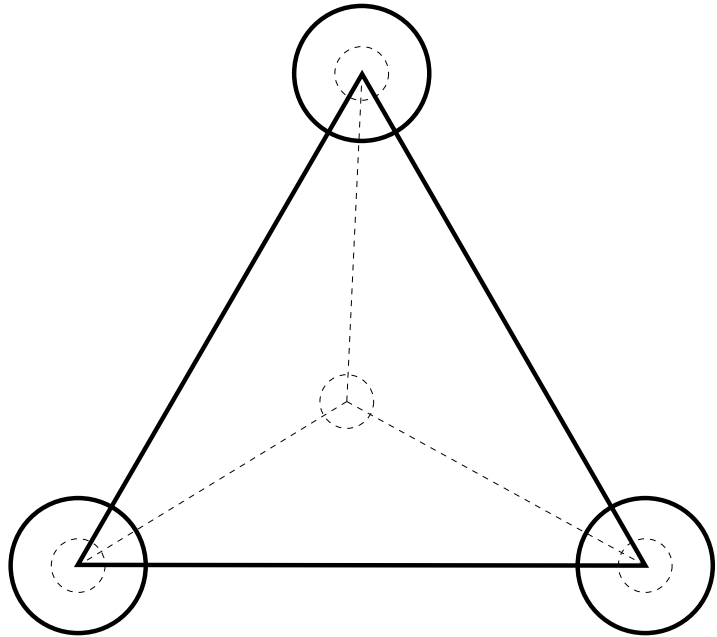

The nets of balls $\mathcal{B}_i$ and $\mathcal{B}_{i+1}$

**Figure 4.1.**

**Proof.** Each ball $\Delta_i \in \mathcal{B}_i$, contains a ball $\Delta_{i+1} \in \mathcal{B}_{i+1}$. Thus $(\mathcal{NS}_d^Q, \mathcal{B}_i) \hookrightarrow (\mathcal{NS}_d^Q, \mathcal{B}_{i+1})$. Moreover if $f \in \mathcal{NS}_d^Q$, then there is some ball $B_1 \subset \Gamma$ such that every point in $B_1$ has $k < Q$ pre-images. $B_1$ contains some ball $\Delta_i \in \mathcal{B}_j$ for $j \gg 0$. Let $x \in \Delta_i$ be a regular value of $f$. Then $f$ is a diffeomorphism restricted to some small ball $B_{y_i}$ surrounding each $y_i \in f^{-1}(x)$. By choosing these balls $B_{y_i}$ small enough we can arrange that they are disjoint. Then for some large $i_\infty$ there is a $\Delta_\infty \in \mathcal{B}_{i_\infty}$ such that

$$\Delta_\infty \subset \cap f(B_{y_i})$$

,and $f \in (\mathcal{NS}_d^Q, \mathcal{B}_{i_\infty})$. $\square$

**Definition 4.16.** *(Bounded Surjectivity Covering) Let $\Delta_k \in \mathcal{N}_i$ Let $U_k \subset (\mathcal{NS}_d^{|d|+1}, \mathcal{N}_i)$ be the maps such that $f^{-1}(\Delta_k)$ consists of $d$ disjoint discs $D_i$, and $f\colon D_i \to \Delta_k$ is a diffeomorphism, and $f$ maps $\Sigma \setminus \bigcup D_i \to \Gamma \setminus \Delta_k$. Then $\bigcup_{k \in S} U_k = (\mathcal{NS}_d^{|d|+1}, \mathcal{N}_i)$.*

**Remark 4.17.** It is straightforward to define analogous coverings for other $Q \neq |d| + 1$, however we will not consider those cases in this paper.

## 4.4 Domain $S^2$: Proofs of Theorems 4.1 and 4.3)

### 4.4.1 Theorem 4.1: $\mathcal{M}_0^\sigma(S^2, \Gamma)$ is homotopy equivalent to $\Gamma$ when $K(\sigma) \leq 1$

Denote the constant maps from $S^2 \to \Gamma$ by $\mathcal{CM}$. Then $\mathcal{CM}$ is homeomorphic to $\Gamma$, and $\mathcal{CM}) \subset (\mathcal{NS}_0^2 T_i)$. We will prove Theorem 4.1 by showing that $\eta\colon \mathcal{CM} \hookrightarrow \mathcal{NS}_0^2(S^2, \Gamma)$ is a homotopy equivalence.

**Proof.** By the Discrete Approximation Lemma (4.15), it is sufficient to show that $f\colon \mathcal{CM}(S^2, \Gamma) \hookrightarrow (\mathcal{NS}_0^2(S^2, \Gamma), \mathcal{B}_i)$ is a homotopy equivalence. We consider the Bounded Surjectivity Cover $U_{\Delta_j}$ of $(\mathcal{NS}_0^2, \mathcal{B}_i)$ provided by Definition 4.16. The sets $V_{\Delta_j} := f^{-1}(U_{\Delta_j})$ then cover $\mathcal{CM}$. I claim that

$$f\colon \bigcap_{j \in S} V_{\Delta_j} \to \bigcap_{j \in S} U_{\Delta_j}$$

is a homotopy equivalence for all indexing sets $S \subset 1..m$. When combined with the Homotopy Decomposition lemma (proposition 4.12) this will show that $f$ must be a homotopy equivalence.

$\bigcap_{j \in S} U_{\Delta_j}$ consists of the maps $S^2 \to \Gamma$ which miss $\bigcup_{j \in S} \Delta_j$. $\bigcap_{j \in S} V_{\Delta_j}$ consists of the constant maps $S^2 \to \Gamma$ which miss $\bigcup_{j \in S} \Delta_j$. Denote $\Gamma, \setminus \bigcup_{j \in S} \Delta_j$ by $\Gamma_\Delta$.

Fix a point $x$ in the domain $S^2$. Then we have a fibration, induced by evaluating each map at $x$:

$$\begin{array}{c} \mathcal{M}_*(S^2, \Gamma_\Delta) \\ \downarrow \\ \bigcap_{j \in S} U_{\Delta_j} \\ \downarrow \pi_x \\ \Gamma_\Delta \end{array}$$

where $\mathcal{M}_*(S^2, \Gamma_\Delta)$ denotes the based maps from $S^2$ to $\Gamma_\Delta$, and $\pi(\gamma) = \gamma(x)$. Evaluation at $x$ also induces a fibration of $\bigcap_{j \in S} V_{\Delta_j}$ over $\bigcap_{j \in S} V_{\Delta_j}$:

$$\begin{array}{c} pt \\ \downarrow \\ \bigcap_{j \in S} V_{\Delta_j} \\ \downarrow \pi_x \\ \Gamma_\Delta \end{array}$$

The inclusion $f$ then induces a morphism of these fibrations:

$$\begin{array}{ccc} f\colon pt & \hookrightarrow & \mathcal{M}_*(S^2, \Gamma_\Delta) \\ \downarrow & & \downarrow \\ f\colon \bigcap_{j \in S} V_{\Delta_j} & \hookrightarrow & \bigcap_{j \in S} U_{\Delta_j} \\ \downarrow \rho_x & & \downarrow \pi_x \\ \Gamma_\Delta & id & \Gamma_\Delta \end{array}$$

$\Gamma_\Delta$ is a surface with nontrivial boundary, it is thus a $K(\pi, 1)$. For any space $X$,

$$\pi_l(\mathcal{M}_*(S^2, X)) \cong \pi_{l+2}(X)$$

Thus as $\Gamma_\Delta$ is a $K(\pi, 1)$, $\mathcal{M}_*(S^2, \Gamma_\Delta)$ is weakly contractible.

So both $\rho_x$ and $\pi_x$ are (weak) homotopy equivalences. Thus

$$f\colon \bigcap_{j \in S} V_{\Delta_j} \hookrightarrow \bigcap_{j \in S} U_{\Delta_j}$$

is also a homotopy equivalence for every indexing set $S$. $f\colon \mathcal{CM}(\mathcal{S}^2, \Gamma) \to (\mathcal{NS}_0^2(S^2, \Gamma), \mathcal{B}_i)$ itself is thus a homotopy equivalence by the Homotopy Decomposition Lemma. $\square$

### 4.4.2 Theorem 4.3: $\mathcal{M}_\sigma^{\pm 1}(S^2, S^2)$ is homotopy equivalent to $SO(3)$

**Proof.** In this section we prove Theorem 4.3:

**Theorem.** *Let $(\Sigma = S^2, \sigma_\Sigma)$ and $(\Gamma = S^2, \sigma_\Gamma)$ be symplectic surfaces. $\mathcal{M}_d^\sigma(S^2, S^2)$ is homotopy equivalent to $SO(3)$ for:*

1. $K(\sigma) \in (0, 1]$, $d = 1$
2. $K(\sigma) \in (1, 2]$, $d = 1$

As $S^2$ admits an orientation reversing involution $\gamma$, it is sufficient to show the first statement. The second will follow by applying the identity of Corollary 2.8. $\mathcal{M}_d^\sigma(S^2, S^2)$ is homotopy equivalent to $\mathcal{NS}_1^3(S^2, S^2)$. $\text{Diff}(S^2)$ is homotopy equivalent to $SO(3)$ [Sma59]. Thus by the Discrete Approximation Lemma (4.15), to prove Theorem 4.3 it is sufficient to show that

$$\eta\colon \text{Diff}(S^2) \hookrightarrow (\mathcal{NS}_1^3, \mathcal{B}_i)$$

is a homotopy equivalence for each $i \in \mathbb{N}$. Denote the Bounded Surjectivity Covering (Definition 4.16) by $U_j$. We aim to prove this Lemma, and thus our theorem, by applying the Homotopy Decomposition Lemma to the $\eta$ and the cover given by $U_j$. Thus we must show that

$$\eta|\bigcap_{j \in S} V_j \to \bigcap_{j \in S} U_j$$

is a homotopy equivalence, where $\{V_j\}$ denotes the cover of $\text{Diff}(S^2)$ given by $\{\eta^{-1}(U_j)\}$. As $U_j \supset \text{Diff}(S^2)$, this cover is trivial: each $V_j$, and thus all of their mutual intersections, consists of all of $\text{Diff}(S^2)$.

Denote $\bigcup_{j \in S} \Delta_j$ by $\Delta$. $\bigcap_{j \in S} U_j$ consists of the maps $S^2 \to S^2$ which are diffeomorphisms on $f^{-1}(\Delta)$, and which map $S^2 \setminus f^{-1}(\Delta)$ to $S^2 \setminus \Delta$. To remind us of its content we will denote $\bigcap_{j \in S} U_j$ by $NS_\Delta$.

Denote $S^2 \setminus \Delta$ by $A_l$ where

$$A_l = D^2 - l \text{ discs}$$

**Definition 4.18.** *Denote by $\mathcal{E}_\Delta$ the space of orientation preserving embeddings $g\colon \Delta \to S^2$.*

**Definition 4.19.** *We denote the degree 1 self maps of $A_l$ which are the identity near the boundary by $(A_l, A_l)^\circ$. Denote the diffeomorphism of $A_i$ which are the identity near the boundary by $\text{Diff}(A_l)^\circ$.*

$\eta$ induces a morphism of fibrations:

$$\begin{array}{ccc}
\eta\colon \text{Diff}(A_l)^\circ & \hookrightarrow & (A_l, A_l)^\circ \\
\downarrow & & \downarrow \\
\eta\colon \text{Diff}(S^2) & \hookrightarrow & NS_\Delta \\
\downarrow \pi & & \downarrow \pi \\
\mathcal{E}_\Delta & id & \mathcal{E}_\Delta
\end{array}$$

Here each projection $\pi$ takes $f \to f^{-1}|_\Delta$. In Lemma 4.20 below we show that

$$\eta\colon \text{Diff}(A_l)^\circ \hookrightarrow (A_l, A_l)^\circ$$

is a homotopy equivalence. Thus, as $\mathcal{E}_\Delta$ is connected:

$$\eta\colon \text{Diff}(S^2) \hookrightarrow NS_\Delta$$

is also a homotopy equivalence by the 5-lemma. □

**Lemma 4.20.** $\eta\colon \text{Diff}(A_l)^\circ \hookrightarrow (A_l, A_l)^\circ$ *is a homotopy equivalence.*

**Proof.** By Theorem 3.1 in [Edm79], $\eta$ induces a bijection on connected components; every $f \in (A, A)_1^\circ$ is homtopic rel $\partial A_1$ to a diffeomorphism in $\text{Diff}(A_l)^\circ$. $\text{Diff}(A_l)^\circ$ consists of contractible components. Thus it is sufficient to show that $(A_l, A_l)^\circ$ also consists of contractible components.

Choose a disjoint set of $l$ arcs $c_i \subset A_l$, such that $A_l \backslash \bigcup c_i$ consists of 2 disjoint discs. Restricting maps to the curves $\bigcup c_i$ then yields a fibration:

$$(A_l, A_l)_2^\circ$$
$$\downarrow$$
$$(A_l, A_l)^\circ$$
$$\downarrow \pi_A$$
$$(A_l, A_l)_1^\circ$$

Where $\pi_A$ sends a map $f$ to its restriction $f|_{\bigcup c_i}$.

$(A, A)_1^\circ$ consists of a direct product of $l$ path spaces in $A_l$. Each path space is homotopy equivalent to the space of based loops in $A_l$. This has contractible components since $A_l$ is a $K(\pi, 1)$. Thus $(A_l, A_l)_1^\circ$ also has contractible components. $(A_l, A_l)_2^\circ$-the extensions of degree one maps, defined on the boundary $\partial A \cup \bigcup_{1..n} c_i$ to all of $A$ - is contractible. We can identify $\pi_i((A_l, A_l)_2^\circ)$ with $\pi_{i+2}(A_l) \oplus \pi_{i+2}(A_l)$. As $A_l$ is a $K(\pi, 1)$ these groups all vanish. Thus $(A_l, A_l)_2^\circ$ is contractible, as are the connected components of $(A_l, A_l)^\circ$. □

## 4.5 Domain $T^2$: (Theorem 4.7) The inclusion $j\colon \Omega_\Gamma \hookrightarrow \mathcal{NS}_0^Q(T^2, \Gamma)$ is a homotopy equivalence for $0 < Q \leq 2$

By Corollary 2.6 is sufficient to consider the case when $Q = 1$.

### 4.5.1 Reduction to based maps

We note that both $\Omega_\Gamma$ and $\mathcal{NS}_0^1(T^2, \Gamma)$ fiber over $\Gamma$ by $f \to f(x)$. We denote these maps by $\pi_\Omega$ and $\pi_N$ respectively. Both are fibrations. The inclusion $j$ then induces a morphism of fibrations:

$$\begin{array}{ccc} (\Omega_\Gamma, y) & \xrightarrow{j_{\text{fb}}} & (\mathcal{NS}_0^1, y) \\ \downarrow & & \downarrow \\ \Omega_\Gamma & \xrightarrow{j} & \mathcal{NS}_0^1(T^2, \Gamma) \\ \downarrow \pi_\Omega & & \downarrow \pi_N \\ \Gamma & \xrightarrow{\text{id}} & \Gamma \end{array}$$

$\Omega_\Gamma$ and $\mathcal{NS}_0^1(T^2, \Gamma)$ are each connected. Thus if we show that the inclusion of the fibers

$$j_{\text{fb}}\colon (\Omega_\Gamma, y) \hookrightarrow (\mathcal{NS}_0^1, y)$$

is a homotopy equivalence, thus inclusion of total spaces

$$j\colon \Omega_\Gamma \hookrightarrow \mathcal{NS}_0^1(T^2, \Gamma)$$

would also be a homotopy equivalence by the 5-lemma. We will prove that $j_{\text{fb}}$ is a homotopy equivalence in the next subsection.

### 4.5.2 Setting up the decomposition

As with the previous theorems make a "discrete approximation" of $\Omega_y$

**Definition 4.21.** Let $\Omega_y^i = \bigl(\mathcal{NS}_0^1(\mathcal{B}_i), y\bigr) \cap (\Omega_\Gamma, y)$

We note that just as the direct system:

$$\hookrightarrow \mathcal{NS}_0^1(\mathcal{B}_i) \hookrightarrow \mathcal{NS}_0^1(\mathcal{B}_{i+1}) ... \hookrightarrow (\mathcal{NS}_0^1, y)$$

(we have suppressed the $\Gamma$ and $T^2$ from our notation) converges to $\mathcal{NS}_0^1$, the direct system:

$$... \hookrightarrow \Omega_y^i \hookrightarrow \Omega_y^{i+1} ... \hookrightarrow (\Omega_\Gamma, y)$$

converges to $\Omega_y$. Thus to prove the lemma it is sufficient to show that the inclusion

$$j\colon \Omega_y^i \hookrightarrow NS_0^1(\mathcal{B}_i)$$

is a homotopy equivalence for each $i$. We apply the Homotopy Decomposition Lemma to this task:

Let $\Delta_{s \in S}$ denote the balls in $\mathcal{NB}^i$. Denote by $U_s$ the maps in $M(T^2, \Gamma)$ which miss the ball $\Delta_s$. This is the Bounded Surjectivity Cover of $NS_0^1(\mathcal{B}_i)$.

$$NS_\Delta = \bigcap_{k \in K} U_k$$

are the $C^1$-maps of $T^2$ into $\Gamma_\Delta = \Gamma \backslash \bigcup_{k \in K} \Delta_k$ which send $NS_\Delta$ $x$ to $y$. On the other hand,

$$\Omega_\Delta = \bigcap_{k \in K} j^{-1}(U_k)$$

is the maps $f \in \Omega_T$ such that $im(f) \subset S_\Delta^2$.

What we must show is that the map:

$$j: \Omega_\Delta \to NS_\Delta$$

is a homotopy equivalence for each indexing set $K \subset S$. This follows from Lemma 4.11 as $\Gamma_\Delta$ is not $S^2$.

### 4.5.3 Proof of Lemma 4.11: If $\Theta \neq S^2$, then the inclusion $j: (\Omega_\Theta, y) \hookrightarrow (\mathcal{M}^0(T^2, \Theta), y)$ is a homotopy equivalence.

**Proof.** We begin by introducing some briefer notation. Let $\mathcal{M}_y^0 = (\mathcal{M}^0(T^2, \Theta), y)$ and let $\Omega_y = (\Omega_\Theta, y)$.

To see this we note that any non spherical surface $\Theta$ is a $K(\pi, 1)$. Thus the based maps $\mathcal{M}_y^0$ consists of contractible components indexed by the possible homomorphisms on the fundamental group:

$$\pi_1(T^2) \hookrightarrow \pi_1(\Gamma)$$

We now characterize these components.

**Definition 4.22.** *Denote by $g \in \mathfrak{G}$ the group homomorphisms*

$$g: \pi_1(T^2) \to \pi_1(\Theta)$$

*which are induced by maps $f \in \mathcal{M}_y^0$.*

**Lemma 4.23.** *Let $\gamma_1$, $\gamma_2$ be as described in Definition 4.5. A homomorphism $g: \pi_1(T^2) \to \pi_1(\Theta)$ is in $\mathfrak{G}$ only if when $g(\gamma_1) \neq 1$, $g(\gamma_2)^a = g(\gamma_1)^b$ for some pair e integers $(a, b) \in \hat{\mathbb{Q}}$. Moreover any mapping $\{\gamma_1, \gamma_2\} \to \pi_1(\Theta)$ satisfying this identity determines a unique $g \in \mathfrak{G}$.*

**Proof.** We begin with the case that $\Theta \neq T^2$. Let $f \in \mathcal{M}_y^0$. Then as $\gamma_1$ and $\gamma_2$ generate $\pi_1(T^2)$ we may describe the homomorphism $f_*$ on $\pi_1$ by determining the values of $f_*(\gamma_1)$ and $f_*(\gamma_2)$. If $f_*(\gamma_1) = \text{id}$ there is no restriction on $f_*(\gamma_2)$. However if $f_*(\gamma_1) \neq \text{id}$ $f_*(\gamma_2)$ is restricted: $\pi_1(T^2)$ is Abelian and thus $[f_*(\gamma_1), f_*(\gamma_2)] = 0$. However two elements $\alpha, \beta \in \pi_1(\Theta)$ commute if and only if there is a pair of relatively prime integers $a, b$ such that $\alpha^a = \beta^b$. Thus $f_*(\gamma_2)^a = f_*(\gamma_1)^b$.

On the other hand, if $\Theta = T^2$ we note that by examining the action on the cohomology rings we can see that any degree 0 self map of the Torus cannot be injective on Homology. The claim follows. $\square$

**Definition 4.24.** *Let $g \in \mathfrak{G}$. Denote by $(\mathcal{M}_y^0, g)$ the $f \in \mathcal{M}_y^0$ such that $f_* = g$. Let $(\Omega_y, g) = (\mathcal{M}_y^0, g) \cap \Omega_y$.*

We note that $(\Omega_\Theta, g)$ is nonempty for each $g \in \mathfrak{G}$. Moreover, as $\Theta$ is a $K(\pi, 1)$ $(\mathcal{M}_y^0, g)$ consists of contractible components. I claim that for each $g \in \mathfrak{G}$, $(\Omega_\Theta, g)$ is nonempty and contractible and thus the inclusion $j_g: (\Omega_\Theta, g) \hookrightarrow (\mathcal{M}_y^0, g)$ is a homotopy equivalence. We denote then the based maps of $(S^1, z)$ to $(\Theta, y)$ by $\mathcal{M}_*(S^1, \Theta)$, and consider two cases:

**Case: $g(\gamma_1) \neq 1$ or $g(\gamma_2) \neq 1$**

In this case there is a unique pair of integers $(a, b) \in \hat{\mathbb{Q}}$ such that $f_*(\gamma_2)^a = f_*(\gamma_1)^b$. Each $f \in (\Omega_\Theta, g)$ factors as $f = h \circ \phi_{a,b}$. $h \in \mathcal{M}_*(S^1, \Theta)$ is any map lying in a homotopy class determined by $g$. Thus $(\Omega_\Theta, g)$ is homeomorphic to a connected component of $\mathcal{M}_*(S^1, \Theta)$ and so contractible as $\Theta$ is a $K(\pi, 1)$.

**Case: $g(\gamma_1) = 1, g(\gamma_2) = 1$**

Each $f$ factors as $f = h \circ \phi_{a,b}$ where $h \in \mathcal{M}_*(S^1, \Theta)$ is a contractible loop. Here $a$ and $b$ may be any pair of relatively prime integers. As $\Theta$ is a $K(\pi_1)$ there is a deformation retract $\lambda_t$ of the contractible loops to the constant loop. Then we can define $\lambda_t(h) \circ \phi_{a,b}$.

□

## 5 Applications to symplectic geometry

By [LM96], symplectic structures on $S^2 \times S^2$ are classified up to diffeomorphism by their cohomology class. Let $(M, \omega_b)$ denote $S^2 \times S^2$ endowed with a symplectic structure $\omega$ such that $K(\omega) = b \geq 1$. We denote the symplectomorphism group of $(M, \omega_b)$ by $\mathcal{S}_b$.

**Definition 5.1.** *Denote by:*

1. $\mathfrak{F}$ *the space of symplectic fibrations of $S^2 \times S^2$ by 2-spheres in the homology class $[\text{pt} \times S^2]$.*

2. $\mathfrak{F}_k$ *the space of pairs $(F, S_k)$ where $F \in \mathfrak{F}$, and $S_k$ is a symplectic section of $F$ in the homology class $k[\text{pt} \times S^2] + [S^2 \times \text{pt}]$*

3. $\mathfrak{F}_{k,-k}$ *the space of triples $(F, S_k, S_{-k})$ where $(F, S_k) \in \mathfrak{F}_k$, and $S_{-k}$ is a symplectic section of $F$ in the homology class $-k[\text{pt} \times S^2] + [S^2 \times \text{pt}]$ which is disjoint from $S_k$.*

4. $\mathcal{E}_k^b$ *the space of symplectic embeddings $\Sigma \to (M, \omega_b)$ which lie in the homology class $k[\text{pt} \times S^2] + [S^2 \times \text{pt}]$*

5. $\mathcal{E}_{k,-k}^b$ *the space of pairs of disjoint embedded symplectic curves $(S_k, S_{-k})$ in $(M, \omega_b)$ such that $S_{\pm k} \in \mathcal{E}_{\pm k}^b$.*

In this section we combine the identities of Corollary 2.8 with work of Abreu and McDuff [AM00] to show:

**Theorem 5.2.** *Let $(M, \omega_b)$ denote $S^2 \times S^2$ endowed with a symplectic structure $\omega$ such that $K(\omega) = b > k \geq 1$ then either the forgetful map $\pi_k: \mathfrak{F}_k \to \mathfrak{F}$ is not a fibration or $\pi_{-k}: \mathfrak{F}_{-k} \to \mathfrak{F}$ is not a fibration.*

**Proof.**

**Idea:** Corollary 2.8 provides certain identities between different spaces of maps with symplectic graphs. We will show that these identities cannot be shared by spaces of symplectic embeddings (which the converse of Theorem 5.2 would predict), by first proving that symplectic embeddings satisfy different identities and then showing a contradiction by comparing with known results on the homotopy type of the symplectomorphism groups. These arguments will be only sketched as the details are very similar to those in [Cof].

### 5.0.4 Almost complex structures and holomorphic fibrations

**Definition 5.3.** *Denote:*

1. *the space of tamed almost complex structures on $(M, \omega_b)$ by $\mathcal{J}$.*

2. *the pairs $(J, S_{\pm k})$ where $S_k \in \mathcal{E}_k^b$, and $J \in \mathcal{J}$ makes $S_{\pm k}$ holomorphic by $\mathcal{J} \pm_k$.*

3. *the triples $(J, S_k, S_{-k})$ where $(S_k, S_{-k}) \in \mathcal{E}_{k,-k}^b$, and $J \in \mathcal{J}$ makes the curves $S_k$ and $S_{-k}$ are holomorphic by $\mathcal{J}_{k,-k}$.*

***Key relationship between almost complex structures and symplectic fibrations:* (Gromov,McDuff-[McD01])** *For each tamed almost complex structure $J$ there is a unique holomorphic fibration by spheres in class $[\text{pt} \times S^2]$.*

The "Key relationship" provides maps: $f\colon \mathcal{J} \to \mathfrak{F}$, $f_k\colon \mathcal{J}_k \to \mathfrak{F}_k$, $f_{k,-k}\colon \mathcal{J}_{k,-k} \to \mathfrak{F}_{k,-k}$, each of which is a fibration with contractible fiber. Further there are also maps $e_k\colon \mathcal{J}_k \to \mathcal{E}_k^b$, and $e_{k,-k}\colon \mathcal{J}_{k,-k} \to \mathcal{E}_{k,-k}^b$. These are also fibrations with contractible fiber. Thus $\mathcal{E}_k^b$ is homotopy equivalent to $\mathfrak{F}_k$, and $\mathcal{E}_{k,-k}^b$ is homotopy equivalent to $\mathfrak{F}_{k,-k}$.

### 5.0.5 Assume $\pi_{\pm k}$ are fibrations; then $\mathcal{E}_{-k}^b \simeq \mathcal{E}_k^{b+2k}$ (will yield contradiction)

By the above $\mathfrak{F}$ is homotopy equivalent to $\mathcal{J}$, and thus contractible. Thus, if $\pi_k$ is a fibration, $\mathfrak{F}_k$ is homotopy equivalent to the symplectic sections of any particular member of $\mathfrak{F}$. In particular we may choose our symplectic fibration to be the trivial one, and we get that $\mathfrak{F}_k$ is homotopy equivalent to $\mathcal{M}_k^{\omega_b}$, and we have:

$$\mathcal{E}_k^b \simeq \mathfrak{F}_k \simeq \mathcal{M}_k^{\omega_b}$$

Similarly $\Sigma_{-k}$ is homotopy equivalent to $\mathcal{M}_b^{-k}$. By our identity (Corollary 2.8) $\mathcal{M}_{-k}^{\omega_b} \approx \mathcal{M}_k^{\omega_b+2k}$, and thus:

$$\mathcal{E}_k^b \simeq \mathfrak{F}_{-k} \simeq \mathcal{M}_{-k}^{\omega_b} \simeq \mathcal{M}_k^{\omega_b+2k} \simeq \mathfrak{F}_k \simeq \mathcal{E}_k^{b+2k}$$

and $\mathcal{E}_k^b \simeq \mathcal{E}_k^{b+2k}$.

### 5.0.6 Identities between spaces of symplectic embeddings

**Proposition 5.4.** *The spaces $\mathcal{E}_k^b$ and $\mathcal{E}_{-k}^b$ are homotopy equivalent.*

**Proof.** We will show that each is homotopy equivalent to $\mathcal{E}_{k,-k}^b$.

$\mathcal{E}_{-k}^b \simeq \mathcal{E}_{k,-k}^b$ We now show that the projection $\mathcal{J}_{k,-k} \to \mathcal{J}_{-k}$ is a homeomorphism. Let $(J, S_{-k}) \in \mathcal{J}_{-k}$, then I claim that there is a unique $J$-holomorphic curve $J$ Fix $k+1$ distinct points $x_i$ on $S_{-k}$. There is a unique $J$-holomorphic curve $F_i$ in class $[F]$ which passes through $x_i$. As both $S$ and the $F_i$ are $J$-holomorphic they must intersect positively. Thus $\Theta$ meets each $F_i$ in precisely one point $\sigma_i$. As $\Theta$ misses $Z_\infty$, $\sigma_i \in F_i - x_i \simeq D^2$. for any $k+1$-tuple in $\prod_{i=1..k+1} F_i - x_i$ there is a unique such curve $\Theta$:

I claim that $\Theta$ is always smooth and irreducible. For as the set of generic almost complex structures is dense one can always approximate $J$ by a sequence of complex structures $J_n$ so that the $J_n$ holomorphic curve through these $q$ points $\Theta_n$ is smooth. The sequence of curves $\Theta_n$ then converges to $\Theta$, and $\Theta$ is thus controlled by Gromov compactness. It consists of a union of $J$-holomorphic spheres, which meet in points. The need to:

1. Intersect the curves in class $[F]$ positively. (Curves in class $[F]$ exist for every $J$ tamed by $\omega$)
2. Intersect $S_{-k}$ positively.

eliminate all such nodal curves, save those of the form:

$$Z_\infty \cup \bigcup_{i=1}^k F_i$$

where the $F_i$ are (possibly repeated) spheres in class $F$. However curves of this last form are eliminated as well. As each point lies off $Z_\infty$ and in a distinct $J$-holomorphic fiber of $F$, the $k$ fiber components $F_i$ cannot pass through all $k+1$ points.

$\mathcal{E}_k^b \simeq \mathcal{E}_{k,-k}^b$ For each tamed almost complex structure $J$ there is a unique holomorphic fibration by spheres in class $[\text{pt} \times S^2]$. This defines maps: $\pi_k\colon \mathcal{J}_k \to \mathfrak{F}_k$, $\pi_{k,-k}\colon \mathcal{J}_{k,-k} \to \mathfrak{F}_{k,-k}$, each of which is a fibration with contractible fiber. Further there are also maps $p_k\colon \mathcal{J}_k \to \mathcal{E}_k^b$, and $p_{k,-k}\colon \mathcal{J}_{k,-k} \to \mathcal{E}_{k,-k}^b$. These are also fibrations with contractible fiber. Thus $\mathcal{E}_k^b$ is homotopy equivalent to $\mathfrak{F}_k$, and $\mathcal{E}_{k,-k}^b$ is homotopy equivalent to $\mathfrak{F}_{k,-k}$.

We will now construct large open neighborhoods of the symplectomorphisms $\mathcal{S}_b$ within the diffeomorphism group. These neighborhood will have the same homotopy type as $\mathcal{S}_b$, but will be far easier to work with. In particular, it will be much easier to understand the "action" of this neighborhood on various objects.

Denote the product fibration of $S^2 \times S^2$ by spheres in class $[\text{pt} \times S^2]$ by $F$. Fix symplectic sections $Z_k$ and $Z_{-k}$ such that the the triple $(F, Z_k, Z_{-k}) \in \mathfrak{F}_{k,-k}$.

**Definition 5.5.** *Denote by $\mathcal{D}_b^{k,-k}$ the diffeomorphisms of $S^2 \times S^2$ which preserve $H^2$ and such that the triple $(g(F), g(Z_k), g(Z_{-k})) \in \mathfrak{F}_{k,-k}$.*

**Proposition 5.6.** *The inclusions $i_+ \colon \mathcal{S}_b \hookrightarrow \mathcal{D}_b^k$ and $i_\pm \colon \mathcal{S}_b \hookrightarrow \mathcal{D}_b^{k,-k}$ are weak deformation retracts.*

The proof of each retraction proceeds in two steps. Denote by $\mathcal{P}_k$ the symplectic forms in class $[\omega_b]$ which are positive on $F$ and $Z_k$. Denote by $\mathcal{P}_{k,-k}$ those in $\mathcal{P}_k$ which are also positive on $Z_{-k}$. Both $P_k$ and $\mathcal{P}_{k,-k}$ are contractible. One shows this by inflation along $Z_k$. Then one gets each retraction by applying Mosers Lemma. For details see Proposition 4.4 in [Cof]. $i_\pm$ is precisely the retraction treated there. $i_+$ can be treated with the same argument.

There are natural orbit maps: $\mathcal{D}_b^k \to \mathfrak{F}_k$, and $\mathcal{D}_b^{k,-k} \to \mathfrak{F}_{k,-k}$. Each of these are surjective fibrations. Denote their fibers by $St_{k,-k}$ and $St_k$ respectively. Each of these are homotopy equivalent to $S^1$. The argument is exactly analogous to those of Propositions 6.1 and 6.2 in [Cof]. The forgetful map $\mathfrak{F}_{k,-k} \to \mathfrak{F}_k$ gives us the following morphism of fibrations:

$$\begin{array}{ccc} \mathcal{S}t_{k,-k} & \stackrel{i_1}{\hookrightarrow} & St_k \\ \downarrow & & \downarrow \\ \mathcal{D}_b^{k,-k} & \stackrel{\text{id}}{\to} & \mathcal{D}_b^k \\ \phi_k \downarrow & & \downarrow \phi_{k,-k} \\ \mathcal{E}_{k,-k}^b \approxeq \mathfrak{F}_{k,-k} & \to & \mathfrak{F}_k \approxeq \mathcal{E}_k^b \end{array}$$

$\mathcal{S}_b$ is connected [AM00]. Thus so are $\mathcal{D}_b^{k,-k}$, $\mathcal{D}_b^k$, $\mathfrak{F}_{k,-k}$ and $\mathfrak{F}_k$. Thus as the maps on the total space and fiber are both weak homotopy equivalences, the 5-lemma implies that the map on the base must also be a homotopy equivalence. □

### 5.0.7 Comparison with rational homotopy of symplectomorphism group yields contradiction

We have the following string of maps comparing the symplectomorphisms $\mathcal{S}_b$ of $(M, \omega_b)$ with those of $(M, \omega_{b+2k})$:

$$\mathcal{S}_b \simeq \mathcal{D}_b^k \underset{\phi_k}{\to} \mathfrak{F}_k \simeq \mathcal{E}_k^b \simeq \mathcal{E}_{-k}^b \simeq \mathcal{E}_k^{b+2k} \simeq \mathfrak{F}_{k,-k} \underset{\phi_{k,-k}}{\leftarrow} \mathcal{D}_{b+2k}^{k,-k} \simeq \mathcal{S}_{b+2k}$$

where each map $\simeq$ is a homotopy equivalence, and $\phi_k$ and $\phi_{k,-k}$ are the orbit maps. As each has fiber homotopy equivalent to $S^1$ they induce isomorphism of $\pi_l$ for $l > 2$. Denote the greater interger $z \le b$ by $\lfloor b \rfloor$. Then by Abreu and McDuff's work on on the rational homotopy type of these symplectomorphism groups *[AM00, McD01] we know that $\pi_{4\lfloor b \rfloor}(\mathcal{S}_b)$ and $\pi_{4\lfloor b \rfloor}(\mathcal{S}_{b+2k})$ are not isomorphic.* Since $b > 1$ this gives our contradiction, and thus completes the proof of Theorem 5.2. □